\documentclass[a4paper]{article}

\usepackage{amsmath,amssymb,graphicx}
\usepackage{subfigure}

\usepackage{lscape}
\usepackage{mathrsfs}

\usepackage{algorithmic}
\usepackage[ruled]{algorithm}

\usepackage{multirow}
\usepackage{verbatim}
\usepackage{enumerate}
\usepackage{color}

\newtheorem{theorem}{Theorem}

\newtheorem{proposition}{Proposition}
\newcommand{\qed}{$\hfill{\Box}$}

\usepackage{longtable}

\title{ M-eigenvalues of The Riemann Curvature Tensor \thanks{To appear in: Communications in Mathematical Sciences.} }

\author{  Hua Xiang$^a$\thanks{
E-mail: hxiang@whu.edu.cn. H. Xiang is supported by the National
Natural Science Foundation of China under grants 11571265, 11471253 and NSFC-RGC 
No.11661161017. }
\quad Liqun Qi$^b$\thanks{
E-mail: maqilq@polyu.edu.hk.  L. Qi is supported by the Hong Kong Research Grant Council
    (Grant No. PolyU  15302114, 15300715, 15301716 and 15300717). }
\quad Yimin Wei$^{c}$\thanks{
E-mail: ymwei@fudan.edu.cn. Y. Wei is supported by the National Natural
Science Foundation of China under grant 11771099 and International Cooperation Project of Shanghai Municipal Science and Technology Commission under grant 16510711200.
}
\\ \\
\small{$^a$ School of Mathematics and Statistics, Wuhan University, Wuhan, 430072, P.R. China}\\
{\small $^b$ Department of Applied Mathematics, The Hong Kong Polytechnic University,
Hong Kong}\\
{\small $^c$ School of Mathematical Sciences and }
\\ {\small Shanghai Key Laboratory of Contemporary Applied Mathematics, }
\\ {\small Fudan University, Shanghai, 200433, P. R. of China}
}

\begin{document}

\date{\today}
\maketitle

\begin{abstract}

The Riemann curvature tensor
is a central mathematical tool in Einstein's theory of general relativity. Its related eigenproblem plays an important role in mathematics and physics.
We extend M-eigenvalues for the elasticity tensor to the Riemann curvature tensor. The definition of M-eigenproblem of the Riemann curvature tensor is introduced from the minimization of an associated function.   The M-eigenvalues of the Riemann curvature tensor always exist and are real.   They are invariants of the Riemann curvature tensor.   The associated function of the Riemann curvature tensor is always positive at a point if and only if the M-eigenvalues of the Riemann curvature tensor are all positive at that point.  We investigate the M-eigenvalues for the simple cases, such as the 2D case, the 3D case, the constant curvature and the Schwarzschild solution, and all the calculated M-eigenvalues are related to the curvature invariants.
\end{abstract}


\section{Introduction}

The eigenproblem of tensor is an very important topic theoretically and practically. In \cite{HanDaiQi09,HuangQi_JCAM18,QDH09,XiangQiWei2017}, the elasticity tensor is investigated, including the strong ellipticity, the positive definiteness, the M-eigenvalues, etc. It is well known that the elasticity tensor is a very important concept in solid mechanics. In this paper, we will consider another counterpart, the Riemannian curvature tensor, which is a basic concept to describe the curved space, and a central mathematical tool in Einstein's general relativity.
In the following, we first review two kinds of eigenproblems associated with the elasticity tensor, and then consider their counterparts corresponding to the Riemannian curvature tensor.

The elasticity tensor $E$ is a fourth-rank tensor. The classical eigenproblem of elasticity tensor reads
\begin{equation}\label{eqn:ElasticityEigentensor}
E_{ijkl} z^{kl} = \zeta ~ z_{ij} ,
\end{equation}
where the eigentensor $z^{ij}$ is symmetric.
For simplicity, when discussing the elasticity tensor, we use the metric in the Kronecker delta, and temporarily omit the difference between the subscripts and the superscripts.

There exists the minor symmetry
$E_{ijkl}=E_{jikl}$ and $E_{ijkl}=E_{ijlk}$, and the major symmetry $E_{ijkl}=E_{klij}$.   The requirement of the symmetry reduces the number of independent elements to 21.
The eigenproblem \eqref{eqn:ElasticityEigentensor}  is closely related to the positive definiteness of $E$,
which has been considered by Lord Kelvin more than 150 years ago,
which also guarantees the uniqueness of solutions in problems of elasticity. 
The elasticity stiffness tensor $E$ must be positive definite, 
which means that the strain energy density or elastic potential satisfies
$$   E_{ijkl} \epsilon^{ij} \epsilon^{kl} > 0, $$
where $\epsilon_{ij}$ is any symmetric strain tensor. It physically means that energy is needed to deform an elastic body from its unloaded  equilibrium position. 
%
The positive definiteness of $E$ requires that all $\zeta$'s are positive.

Qi, Dai and Han \cite{QDH09} in 2009 introduced M-eigenvalues for the elasticity tensor.
The M-eigenvalues $\theta$ of the fourth-order tensor are defined as follows.
\begin{eqnarray} \label{eqn:MeigDef_x&y}
 E_{ijkl} y^j x^k y^l = \theta x_i , \quad
 E_{ijkl} x^i y^j x^k = \theta y_l ,
\end{eqnarray}
under the constraints $\langle x, x \rangle = \langle y, y \rangle = 1$, i.e., $x$ and $y$ are normalized.  Here, $x$ and $y$ are real vectors.   The M-eigenvalues always exist and are real.
Furthermore, they are isotropic invariants of the elasticity tensor.   For the elasticity tensor $E$, since $E_{ijkl} x^i y^j x^k = E_{kjil} x^i y^j x^k$, and $ E_{kjil}= E_{lijk} $, the second equality is equivalent to
$$E_{ijkl} x^j y^k x^l = \theta y_i . $$

For the fourth-order elasticity $E$,
the strong ellipticity is defined by the following function to be positive.
\begin{equation*}
f(x,y) = E_{ijkl} x^i y^j x^k y^l > 0, \ \ \forall  x,y \in \mathbb{R}^3.
\end{equation*}
Such strong ellipticity condition ensures that the governing differential equations for elastostatics problems be completely elliptic. It is an important property in the elasticity theory associated with uniqueness, instability, wave propagation, etc, and has been studied extensively.
It was shown that the strong ellipticity condition holds if and only if all the M-eigenvalues are positive.


These two kinds of eigenproblems relate to the positive definiteness and the strong ellipticity of the material respectively.
The positive definiteness is a less general hypothesis than the strong ellipticity.
The positive definiteness  implies strong ellipticity, while the converse statement is not true.
%
Motivated by the work on the elasticity \cite{QDH09,HanDaiQi09,XiangQiWei2017,QiChenChen18,HuangQi_JCAM18}, in this paper we consider the corresponding tensor eigenproblems for the Riemann tensor.

Let $(M,g)$ be an $n$-dimensional Riemannian manifold. That is, $M$ is the Riemannian manifold equipped with the Riemannian metric $g$.
We consider the curvature tensor of the Levi-Civita connection $\nabla$ of the Riemannian metric $g$.
The curvature tensor for the Levi-Civita connection will be called later the Riemann curvature tensor, or the Riemann tensor.

The curvature $R$ of a Riemannian manifold $M$  corresponds to a mapping $R(X,Y)$  associated  to the pair $(X,Y)$  by
$$R(X,Y) :=  [\nabla_X, \nabla_Y] - \nabla_{[X,Y]} .$$

Let $T_p M$ be the tangent space of $M$ at the point $p$.
The (0,4) type Riemannian curvature tensor is a quadrilinear mapping:
$$ R: ~ T_p M \times T_p M \times T_p M \times T_p M  \rightarrow \mathbb{R} .$$
$$ R(W, Z, X, Y) := \langle W, R(X,Y) Z  \rangle, \quad \forall ~W, X, Y, Z \in T_p M .$$
By the way, the (1,3) type is given by
$(\omega, Z, X, Y) \mapsto \omega( R(X,Y) Z )$, $\forall$ vector fields $ X, Y, Z$ and 1-form $\omega$. Here
$R(X,Y) Z$ or $R(W, Z, X,Y)$ is called the curvature tensor of the Levi-Civita connection.   The notation here is somewhat abused.   The two mappings $R(X,Y$ (also denoted by $R_{XY}$ in some references) and $R(W, Z, X,Y)$ use the same letter ``R''.  As such a usage can be found in the literature \cite{BaezMuniain_Book1994, Boothby_Book86, ChernChenLam_book1999, DoCarmo_book92, Lee_Book97} and will not cause confusion, we keep such a usage.


To work with components, one needs a local coordinate $\{ x^i \}$, a set of corresponding basis vector $\{ \partial_i \}$ and the dual set of basis 1-forms $\{ dx^i \}$.
Suppose that $g_{ij} := g( \partial_i, \partial_j )$ and define Christoffel symbols of the Levi-Civita connection  by the formula
$$\Gamma^i_{jk} 
= \frac{1}{2} g^{ih} \left( \frac{\partial g_{hj}}{\partial x^k} + \frac{\partial g_{hk}}{\partial x^j} - \frac{\partial g_{jk}}{\partial x^h}  \right). $$

%

Using $[\partial_i, \partial_j] = 0$ and $\nabla_{\partial_i} \partial_k = \Gamma^l_{ik} \partial_l$, we can calculate \cite[P.25]{Wolf_book6th}
\begin{eqnarray*}
R( \partial_i, \partial_j ) \partial_k &=& ( [\nabla_{\partial_i} , \nabla_{\partial_j}] - \nabla_{[\partial_i,\partial_j]}) \partial_k  = \nabla_{\partial_i} \nabla_{\partial_j} \partial_k - \nabla_{\partial_j} \nabla_{\partial_i} \partial_k  \\
&=& ( \partial_i \Gamma^l_{jk} - \partial_j \Gamma^l_{ik}   + \Gamma^h_{jk}\Gamma^l_{ih} - \Gamma^h_{ik} \Gamma^l_{jh}) \partial_l  \equiv  R^l_{kij} \partial_l  .
\end{eqnarray*}
That is, $R^l_{kij} = dx^l (R( \partial_i, \partial_j ) \partial_k) $.
For convenience, one can consider the full covariant Riemann curvature tensor  $R_{ijkl} $.
$$ R(\partial_i, \partial_j, \partial_k, \partial_l) = \langle  \partial_i,  R(\partial_k,\partial_l) \partial_j \rangle = g_{ih} R^h_{jkl} \equiv  R_{ijkl}. $$
In terms of the Riemann metric and the coefficients $\Gamma_{ij}^k$ of the Riemannian connection,   we have
$$R_{ijkl}  = \frac{1}{2} \left( g_{il,jk} - g_{ik,jl} + g_{jk,il} - g_{jl,ik} \right)
+g_{hm} ( \Gamma^h_{il} \Gamma^m_{jk} - \Gamma^h_{ik} \Gamma^m_{jl} ). $$

Let $W = w^i \partial_i$, $Z = z^j \partial_j$, $X = x^k \partial_k$, $Y = y^l \partial_l$.
We can verify that
$$R(X,Y) Z = R(x^k \partial_k, y^l \partial_l) (z^j \partial_j) = z^j x^k y^l  R( \partial_k, \partial_l)  \partial_j =  z^j x^k y^l R^h_{jkl} \partial_h, $$
and
$$R(W,Z,X,Y) = w^i z^j x^k y^l R_{ijkl}. $$

The curvature tensor has the following symmetry properties \cite{Boothby_Book86,DoCarmo_book92,Lee_Book97}:
$$R_{ijkl} = -R_{jikl} = -R_{ijlk} = R_{klij}, \quad  R_{ijkl} + R_{iljk} + R_{iklj} =0 .$$
The second identity is called the first (or algebraic) Bianchi identity.
For the 4D case, there are 256 components, but only twenty are independent because of these symmetries \cite{BaezMuniain_Book1994,Dirac_Book75,MTW1973}.

The contraction yields the Ricci tensor $R_{ik}$ and the Ricci scalar $R$ as follows.
$$ R_{ik} = g^{hj} R_{hijk} = R^m_{imk} =  \partial_l \Gamma^l_{ik} - \partial_k \Gamma^l_{il} + \Gamma^l_{ik} \Gamma^h_{lh} - \Gamma^h_{il} \Gamma^l_{hk} ,   $$
$$ R =  R^k_k = g^{ik} R_{ik} = g^{ik} R^m_{imk} = g^{ik} g^{hj} R_{hijk}. $$

In the following of the paper, we present two kinds of eigenproblems in Section 2, just as the eigenproblems for the elasticity tensor. In Section 3, we study four typical cases, calculate the M-eigenvalues and examine their relationship with some well-known invariants.

\section{Eigenproblems of the Riemann tensor}

We consider two kinds of eigenproblems associated with the Riemann tensor.
The first one is well-studied. For the invariant characterizations of a gravitational field, to investigate the algebraic structure of the tensor we consider the eigenproblem \cite{Stephani_Book2003}
\begin{equation}\label{eqn:RCurvatureEigentensor}
R_{ijkl} x^{kl} = \zeta ~ x_{ij} = \zeta ~ g_{im} g_{jn} x^{mn} ,
\end{equation}
where the eigentensor $x^{ij}$ is anti-symmetric.


To express it in a compact matrix form,  we identify a pair of indices $\{ij\}$ of 4D indices with a multi-index that has the range from 1 to 6: 10 $\rightarrow$ 1, 20 $\rightarrow$ 2, 30 $\rightarrow$ 3, 23 $\rightarrow$ 4, 31 $\rightarrow$ 5, 12 $\rightarrow$ 6.
Denoting the basis indices by capital letters and using the symmetries, we can rewrite the above eigenproblem  \eqref{eqn:RCurvatureEigentensor} as follows \cite{Petrov1954}.
\begin{equation}\label{eqn:RCurvatureEig_inMatrix}
R_{AB} x^B = \frac{1}{2} \zeta ~ G_{AB} x^B,
\end{equation}
where both $ ( R_{AB} )$ and $ ( G_{AB} )$ are 6-by-6 matrices, and $ ( R_{AB} )$ can be further  expressed by two symmetric 3-by-3 matrices according to the Einstein field equation in vacuum (see Appendix B).
For this vacuum case, the Riemann tensor is equivalent to the Weyl tensor.
And the well-known Petrov classification reduces to investigate the eigenvalues and the independent eigenvectors and results in Petrov Types I (D), II (N) and III \cite[P.235]{Padmanabhan_Book2013}.


Next, we introduce the M-eigenproblem of the Riemann curvature tensor. Let us consider the associated function
$$ Q(u,v)  \equiv    \langle u, R(u,v) v \rangle = R(u, v, u, v) = R_{ijkl} u^i v^j u^k v^l   . $$

The Riemann curvature tensor is a real tensor.  Vectors $u$ and $v$ are real vectors.   Thus, $Q(u, v)$ is a real continuous function, and we may
consider the following optimization problem
\begin{eqnarray*}
\min_{u,v \in T_p M }  Q(u, v) \qquad
\text{s.t.} \quad \langle u,u \rangle = \langle v,v \rangle = 1.
\end{eqnarray*}
The feasible set of this optimization problem is compact.  Hence, it always has global optimal solutions.  It has only equality constraints.   By optimization theory, its optimal Lagrangian multipliers $\lambda$ and $\mu$ always exist and are real.
The optimality condition reads
\begin{eqnarray*}
R(v,u) v &=& \lambda  u,  \\
R(u,v) u &=& \mu      v,  \\
\langle u,u \rangle &=& 1, \\
\langle v,v \rangle &=& 1.
\end{eqnarray*}
This ensures the existence of such $\lambda$ and $\mu$, and they are real.
It is easy to verify that $\lambda = -R(u,v,u,v) = \mu $.  We may rewrite the M-eigenproblem as seeking the normalized eigenvector pair $(u,v)$ and the eigenvalue $\theta$ in a coordinate-free manner as follows.
$$\widetilde{R} (u,v) \equiv \left( R(u,v) v, R(v,u) u \right)  =   \theta ( u, v ),$$
i.e.,
\begin{align}
R(v,u) u = \theta v ,   ~~ & ~~ R(u,v) v =  \theta u , \label{eqn:Meig_curvature_tensor} \\
\langle u,u \rangle = 1, \quad & \quad
\langle v,v \rangle = 1.
\end{align}
Thus, the M-eigenvalues always exist and are real.
The corresponding component form reads
$$ R^l_{ijk} u^i v^j u^k =  \theta   v^l ,  \quad
   R^l_{ijk} v^i u^j v^k =  \theta   u^l,  $$
or, equivalently,
\begin{eqnarray}\label{eqn:RCurvature_Meig_component}
R_{hijk} u^i v^j u^k =  \theta g_{hl} v^l  = \theta v_h  ,  \quad
R_{hijk} v^i u^j v^k =  \theta g_{hl} u^l  = \theta u_h  ,
\end{eqnarray}
where $g_{ij}u^i u^j = g_{ij}v^i v^j =1$.

Using $\langle u,u \rangle  = \langle v,v \rangle =  1$, we have the M-eigenvalue
\begin{equation}\label{eqn:Meig_theta}
 \theta =    \langle u, R(u,v) v \rangle = R(u, v, u, v) = R_{ijkl} u^i v^j u^k v^l   .
\end{equation}

In some cases we further require that $\langle u,v \rangle =0$. We add this constraint and  present the following modified M-eigenvalue problem.
\begin{eqnarray}
R(v,u) u = \theta v , \quad R(u,v) v =  \theta u ,  \label{eqn:modifiedMeig} \\
\langle u,u \rangle =
\langle v,v \rangle = 1, ~~~
\langle u,v \rangle = 0.  \label{eqn:modifiedMeig_constraint}
\end{eqnarray}
Similarly, the modified M-eigenvalues always exist and are real.

Here the equation (\ref{eqn:RCurvature_Meig_component}) is a tensor equation.  Hence, the M-eigenvalues are invariants of the Riemann curvature tensor.  This is also true for the modified M-eigenvalues.    Also, from the properties of the above optimization problem, the associated function $Q(u, v)$ is always positive at a point if and only if all the M-eigenvalues of the Riemann curvature tensor at that point are all positive.   We see that the M-eigenvalues of the Riemann curvature tensor has all the good properties of the M-eigenvalues of the elasticity tensor.

{\bf Remark 1.}
The optimization problem is just used to introduce the M-eigenvalue problem on Riemannian manifold. The definition should not be limited on a Riemannian manifold, and it can be extended to a pseudo-Riemannian manifold. But for the pseudo-Riemannian manifold, there may be no maximum or minimum of the optimization problem, since the constraint set is not compact any more. In the Lorentz manifold, if $v$ and $u$ are spacelike and timelike orthonormal vectors respectively, i.e.,
$\langle v,v \rangle = 1$ and $\langle u,u \rangle = -1$, then  we have
$\theta = 0$, since the first formula in \eqref{eqn:modifiedMeig} yields $ \theta = \langle v, R(v,u) u \rangle = R(v,u,v,u)$, while the second formula in \eqref{eqn:modifiedMeig} gives $\theta = - \langle u, R(u,v) v \rangle  = - R(u,v,u,v)$, where $R(u,v,u,v)=R(v,u,v,u)$ due to the symmetry property.

{\bf Remark 2.}
Suppose that $u=v$.   Then the M-eigenvalue problem reduces to $R(u,u)u = \theta u$ with $\langle u,u \rangle = 1$. In componentwise form, it can be written as $R_{ijkl} u^j u^k u^l = \theta u_i$, which is a Z-eigenvalue problem \cite{Qi_JSC05}. Since $R_{ijkl} u^j u^k u^l = - R_{ijlk} u^j u^l u^k \equiv 0$, we have $\theta=0$.

Recall that we have two kinds of eigenproblems \eqref{eqn:ElasticityEigentensor} and \eqref{eqn:MeigDef_x&y} for the elasticity tensor.  For the Riemann tensor, we have similar things, i.e., \eqref{eqn:RCurvatureEigentensor} and \eqref{eqn:RCurvature_Meig_component}.
%
In a special case, the M-eigenvalue in \eqref{eqn:RCurvature_Meig_component} is related to the classical eigenvalue in \eqref{eqn:RCurvatureEigentensor} as stated below.
\begin{theorem}
Suppose that   $(\zeta, x)$ is the eigenpair of \eqref{eqn:RCurvatureEigentensor} with $x^{ij} = u^i v^j - v^i u^j$ and $\langle u, v \rangle = 0$. Then $(\theta, u, v)$ is the eigentriple of modified M-eigenproblem  \eqref{eqn:modifiedMeig}-\eqref{eqn:modifiedMeig_constraint}  and $\zeta = 2 \theta$.
\end{theorem}
{\bf Proof}.
Substituting $x^{ij} = u^i v^j - v^i u^j$, the eigenproblem \eqref{eqn:RCurvatureEigentensor} reads
$$ R_{ijkl} (u^k v^l - v^k u^l) = \zeta ( u_i v_j - v_i u_j ) . $$
It is easy to verify that
$$ R_{ijkl} v^j u^k v^l - R_{ijkl} v^j v^k u^l = \zeta ( u_i v_j v^j - v_i u_j v^j ) . $$
Using $R_{ijkl} v^j v^k u^l = R_{ikjl} v^j v^k u^l = - R_{iklj} v^k u^l v^j = - R_{ijkl} v^j u^k v^l $, we have
\begin{equation}\label{eqn:proofThm1-1}
    R_{ijkl} v^j u^k v^l =  \frac{1}{2} \zeta ( u_i - v_i \langle u,v \rangle ) .
\end{equation}
Note that the eigenproblem  \eqref{eqn:RCurvatureEigentensor}  is equivalent to
$R_{ijkl} x^{ij} = \zeta g_{km} g_{ln} x^{mn} $.
Similarly,   we have
\begin{equation}\label{eqn:proofThm1-2}
    R_{ijkl} u^i v^j u^k   = \frac{1}{2}  \zeta ( v_l - u_l \langle u,v \rangle ) .
\end{equation}

Using the orthogonality that $\langle u,v \rangle = 0$, then \eqref{eqn:proofThm1-1} and \eqref{eqn:proofThm1-2} reduce to
$$
R_{ijkl} v^j u^k v^l =  \frac{1}{2} \zeta   u_i  , \quad
R_{ijkl} u^i v^j u^k   = \frac{1}{2}  \zeta  v_l   .
$$
This is just the modified M-eigenproblem \eqref{eqn:modifiedMeig} with the M-eigenvalue $\theta = \frac{1}{2}  \zeta $.

\qed

Besides,
using $\langle u,u \rangle  = \langle v,v \rangle =  1$ and the formulas \eqref{eqn:proofThm1-1} and \eqref{eqn:proofThm1-2}, we have
$$ \frac{\zeta}{2} = \frac{ R_{ijkl} u^i v^j u^k v^l }{ 1 - \langle u,v \rangle^2 } = \frac{ R(u,v,u,v) }{ 1 - \langle u,v \rangle^2 } . $$
The orthogonality that $\langle u,v \rangle = 0$ and the formula \eqref{eqn:Meig_theta} again yield that $\zeta = 2 \theta$.

The sectional curvature is closely related to the M-eigenvalue.
Let $\pi \subset T_p M$ be a 2D subspace of the tangent space $T_p M$ and let $u, v \in \pi$ be two linearly independent vectors (not necessarily orthonormal).
The sectional curvature of $(M,g)$ at $p$ with respect to the 2D plane  $\pi = \text{span} \{u, v\} \subset T_p M$, independent of the choices of basis $\{ u , v \}$, is defined by
\begin{equation}\label{eqn:sectionCurvature}
K(\pi)   = \frac{R(u,v,u,v)}{|u \wedge v |^2} ,
\end{equation}
where
$ |u \wedge v|^2 := \langle u,u \rangle \langle v,v \rangle - \langle u,v \rangle^2 = g_{ij} u^i u^j g_{kl} v^k v^l - (g_{ij} u^i v^j)^2  $ denotes the square of the area of the 2D parallelogram spanned by the pair of vectors $u$ and $v$.

\begin{theorem}
Using the notation above, the eigenvalue of modified M-eigenproblem \eqref{eqn:modifiedMeig}-\eqref{eqn:modifiedMeig_constraint} is the sectional curvature $K(\pi)$.
\end{theorem}
{\bf Proof.}
Let  $(\theta, u, v)$ be the eigentriple of modified M-eigenproblem  \eqref{eqn:modifiedMeig}-\eqref{eqn:modifiedMeig_constraint}.
Taking the inner product with the first equation of \eqref{eqn:modifiedMeig}, we obtain
$\langle v, R(v,u)u \rangle = \theta \langle v, v \rangle$. Using $\langle v, v \rangle = 1$, we have $\theta = R (v, u, v, u)$.
Similarly, from the second equation of \eqref{eqn:modifiedMeig}, we have  $\theta = \langle u, R(u,v) v \rangle = R(u,v,u,v)$.

When $(u,v)$ is chosen as an orthonormal basis for $\pi$, we have $|u \wedge v|^2 =1$, and \eqref{eqn:sectionCurvature} reduces to
$K(\pi) = R ( u, v, u, v ) = R_{ijkl} u^i v^j u^k v^l $, and hence $\theta = K(\pi)$.

\qed

The sectional curvature is essentially the restriction of the Riemann curvature tensor to special set of vectors. 
The knowledge of $K(\pi)$, for all $\pi$, determines the curvature $R$ completely \cite[P.94]{DoCarmo_book92}.
Due to Theorem 2, it is not surprising to see that the M-eigenvalues of the cases examined in the next section relate to the Riemann curvature scalar $R$.

In the following we point out another related thing:  the Jacobi equation, which reads
\begin{equation}\label{eqn:Jacobi}
 \nabla_u \nabla_u v + R(v,u) u = 0,
\end{equation}
where $u$ is the tangent vector and $v$ is a vector field along the geodesic \cite[P.310]{ChernChenLam_book1999}. Each Jacobi field tells us how some family of geodesics behaves \cite{Lee_Book97}.
In component notation, it can be written as  
$$ \frac{D^2 v^\mu}{ d t^2} + R^\mu_{\nu \rho \sigma} u^\nu v^\rho u^\sigma = 0,  $$
where $t$ represents an affine parameter along the geodesics.

Note that the second term of the Jacobi equation   also appears in the M-eigproblem \eqref{eqn:Meig_curvature_tensor}, where it reads
\begin{equation}\label{eqn:Jacobi2nd_Meig1st}
 R^\mu_{\nu \rho \sigma} u^\nu v^\rho u^\sigma =  \theta v^\mu.
\end{equation}
Using \eqref{eqn:Jacobi2nd_Meig1st},  we can express the Jacobi equation as
$$  \ddot{v}^\mu  + \theta v^\mu = 0,  $$
where the dot denotes the ordinary derivative with respect to $t$.
This is a linear system of second-order ODEs for the functions $v^\mu$, and it can be solved by using the proper initial conditions \cite{DoCarmo_book92}.
%

Let us deviate for a while, and switch to the case in the Lorentz manifold, where we define $u$ as the tangent vector to the geodesic and $v$ the geodesic separation, the displacement from fiducial geodesic to nearby geodesic with the same affine parameter. Then \eqref{eqn:Jacobi} is also called as the equation of geodesic deviation \cite[P.219]{MTW1973}.
This equation gives the relative acceleration of free particles \cite{Pirani_APP56}, can server as a definition of the Riemann curvature tensor whose components can be determined
by throwing up clouds of test particles and measuring the relative
accelerations between them \cite{Szekeres_JMP65}.
The  eigenproblem \eqref{eqn:Jacobi2nd_Meig1st} associated with  the second term of the equation of geodesic deviation, which also appears in the M-eigproblem, has explicit physical meaning.
For a rigid body in free fall, the nontrivial eigenvalues of \eqref{eqn:Jacobi2nd_Meig1st} give the principal internal stresses to keep all the parts of the body together in a rigid shape \cite{Padmanabhan_Book2013}.
This may shed some light on the physical meaning of the M-eigenvalues of the Riemann tensor.


\section{Case study}

In this section we calculate several simple concrete cases including the conformally flat case and the Schwarzschild solution to examine what the M-eigenvalues are.

\subsection{The 2D case }

For a 2D case, the Riemann tensor can be expressed by $R_{abcd} = K (g_{ac} g_{bd} - g_{ad}g_{bc}) $, where $K$ is called the Gaussian curvature \cite[P.101]{Stephani_Book2003}.
And there is only one independent component $R_{1212}$.

We can easily check that
$$R_{11} = R^i_{1i1} = R^2_{121} = g^{21}R_{1121} + g^{22}R_{2121}=  g^{22} R_{1212},$$
$$R_{22} = R^i_{2i2} = g^{11} R_{1212}, ~
 R_{12} = R_{21} = R^i_{1i2} = - g^{12} R_{2121} . $$
Then the scalar curvature $R$ can be expressed by
\begin{eqnarray*}
R &=& R_{kj} g^{kj} = g^{11} R_{11} + 2g^{12} R_{12} + g^{22} R_{22}  \\
  &=&  2 R_{1212} (g^{11} g^{22} - g^{12} g^{12}) =  2 R_{1212} \det g^{-1} .
\end{eqnarray*}

The Gaussian curvature, whose magnitude defined by external observer, equals to the scalar curvature, the magnitude defined in terms of internal observer. That is,
$$K =  \frac{R_{1212}}{ \det g} = \frac{R}{2} .$$
This is nothing but Theorema Egregium of Gauss.

\begin{proposition}
For the 2D case, the nonzero M-eigenvalue is $\zeta = K$, where $K$ is the Gaussian curvature.
\end{proposition}
{\bf Proof}.
In components, the M-eigenvalue problem reads
\begin{eqnarray*}
R_{1212} ( y^2 x^1 y^2 - y^2 x^2 y^1 ) &=& \zeta x_1 = \zeta g_{1k} x^k = \zeta ( g_{11} x^1 + g_{12} x^2 ), \\
R_{1212} ( y^1 x^2 y^1 - y^1 x^1 y^2 ) &=& \zeta x_2 = \zeta g_{2k} x^k = \zeta ( g_{21} x^1 + g_{22} x^2 ), \\
R_{1212} ( x^2 y^1 x^2 - x^2 y^2 x^1 ) &=& \zeta y_1 = \zeta g_{1k} y^k = \zeta ( g_{11} y^1 + g_{12} y^2 ), \\
R_{1212} ( x^1 y^2 x^1 - x^1 y^1 x^2 ) &=& \zeta y_2 = \zeta g_{2k} y^k = \zeta ( g_{21} y^1 + g_{22} y^2 ),
\end{eqnarray*}
with the constraints
$$\langle x,x \rangle  = g_{ka} x^a x^k= g_{11} x^1 x^1 + 2g_{12}x^1 x^2 + g_{22} x^2 x^2 = 1,$$
$$\langle y,y \rangle   = g_{ka} y^a y^k= g_{11} y^1 y^1 + 2g_{12}y^1 y^2 + g_{22} y^2 y^2 = 1. $$

Solving this system of polynomial equations, we have $\zeta = 0$ and
$$ \zeta = \frac{ R_{1212} }{ g_{11}g_{22}-g_{12}^2 } = \frac{R}{2} = K. $$
\qed

For the modified M-eigenvalue problem, we add the constraint
$$\langle x,y \rangle   = g_{ij} x^i y^j= g_{11} x^1 y^1 + g_{12}x^1 y^2 + g_{21}x^2 y^1 + g_{22} x^2 y^2 = 0. $$
We then have $\zeta = K$ without the zero eigenvalue.

\subsection{The 3D case }

For the 3D case, there are six independent components in the Riemann tensor, and it can be expressed by Ricci tensor as follows.
\begin{equation} \label{eqn:Rijkl_3D}
 R_{abcd} = R_{ac}g_{bd} - R_{ad}g_{bc} + g_{ac} R_{bd} - g_{ad} R_{bc} - \frac{R}{2}
( g_{ac}g_{bd} - g_{ad}g_{bc} ).  
\end{equation}

Suppose that $(\lambda, x)$ and $(\mu, y)$ are the eigenpairs of the Ricci tensor ($x \neq y$), that is, $R_{ac} x^c = \lambda x_a$ and $R_{ad} y^d  = \mu y_a$ \cite{Hamilton_JDG82,Hamilton_JDG86}. 
The M-eigenvalues of the Riemann tensor are related to the eigenvalues of the Ricci tensor as stated in the following proposition.

\begin{proposition}
For the 3D case, if $\lambda$ and $\mu$ are the eigenvalues of the Ricci tensor associated with eigenvectors $x$ and $y$ respectively, then $(\zeta, x, y)$ is the eigentriple of the modified M-eigenproblem with the M-eigenvalue $\zeta = \lambda + \mu -\frac{R}{2}$.
\end{proposition}
{\bf Proof}.
For the M-eigenvalue problem, we need to calculate
\begin{equation} \label{eqn:RC3D0}
R_{abcd} y^b x^c y^d = \zeta g_{ak} x^k  , ~ (a=1,2,3).
\end{equation}
Substituting the expression of \eqref{eqn:Rijkl_3D} into the left hand side (LHS) of \eqref{eqn:RC3D0}, we have
\begin{eqnarray*}
\text{LHS of \eqref{eqn:RC3D0}}
&=& R_{ac} x^c - \langle x,y \rangle R_{ad} y^d + R_{bd} y^b y^d x_a - R_{bd} y^b x^d y_a - \frac{R}{2} \left( x_a - \langle x,y \rangle y_a \right) \\
&=& R_{ac} x^c - \langle x,y \rangle R_{ad} y^d + y^b R_{bd} (x^c y^d - x^d y^c) g_{ca} - \frac{R}{2} \left( x_a - \langle x,y \rangle y_a \right).
\end{eqnarray*}
Since $R_{ac} x^c = \lambda x_a$,   $R_{ad} y^d = \mu y_a$, and $y^b y_b=1$, the 3rd term equals to
$$ y^b (\mu y_b x^c - \lambda x_b y^c) g_{ca} = \mu x_a - \lambda \langle x,y \rangle y_a .$$

Direct calculation yields that
\begin{eqnarray*}
\text{LHS of \eqref{eqn:RC3D0}}
= (\lambda + \mu) ( x_a -   \langle x,y \rangle y_a ) - \frac{R}{2}  \left( x_a - \langle x,y \rangle y_a \right) .
\end{eqnarray*}
Since $R_{ab}$ is symmetric, if $\lambda \neq \mu$, then $\langle x,y \rangle = 0 $; if $\lambda = \mu$, one can orthogonalize the vectors such that  $\langle x,y \rangle = 0 $. Hence, $\text{LHS of \eqref{eqn:RC3D0}}
= (\lambda + \mu -\frac{R}{2}) x_a $.

The second equation of the M-eigenproblem reads
\begin{equation} \label{eqn:RC3D1}
R_{abcd} x^b y^c x^d = \zeta g_{ak} y^k  , ~ (a=1,2,3).
\end{equation}
We can calculate that
\begin{eqnarray*}
\text{LHS of \eqref{eqn:RC3D1}}
&=& R_{ac} y^c - \langle x,y \rangle R_{ad} x^d + R_{bd} x^b x^d y_a - R_{bd} x^b y^d x_a - \frac{R}{2} \left( y_a - \langle x,y \rangle x_a \right) \\
&=& R_{ac} y^c - \langle x,y \rangle R_{ad} x^d + x^b R_{bd} (y^c x^d - y^d x^c) g_{ca} - \frac{R}{2} \left( y_a - \langle x,y \rangle x_a \right) .
\end{eqnarray*}
Similarly we have that
$\text{LHS of \eqref{eqn:RC3D1}}
= ( \mu + \lambda - \frac{R}{2} ) y_a  $.

Comparing with the right hand side of \eqref{eqn:RC3D0} and \eqref{eqn:RC3D1}, we
 have $\zeta = \lambda + \mu -\frac{R}{2}$.

\qed

\subsection{The constant curvature}

The curvature tensor of a space of constant curvature is expressed in terms of the curvature $\kappa$ and the metric tensor $g_{ij}$  by the formula \cite{Synge_Annals34}
\begin{equation} \label{eqn:ConstCurvatureR}
 R^i_{jkl} = \kappa ( g^i_k g_{jl} - g^i_l g_{jk} ) .
\end{equation}
For $\kappa=0$,  it is an Euclidean space; for $\kappa>0$, it is the sphere of radius $1 / \sqrt{\kappa}$; and for $\kappa<0$, it is a Lobachevskii space.

\begin{proposition}
For the spaces of constant curvature \eqref{eqn:ConstCurvatureR}, $\zeta = \kappa$ is the M-eigenvalue.
\end{proposition}
{\bf Proof}.
Suppose that $x,y$ are the M-eigenvectors. Direct calculation yields that
$$R_{abcd} y^b x^c y^d = \kappa ( g_{ac} g_{bd} y^b x^c y^d - g_{ad} g_{bc} y^b x^c y^d ) = \kappa ( x_a - \langle x,y \rangle y_a ) = \zeta x_a .$$

Similarly, we have $R_{abcd} x^b y^c x^d =  \kappa ( y_a - \langle x,y \rangle x_a ) = \zeta y_a $.

Using $\langle x,y \rangle =0$, then we have the M-eigenvalue $\zeta = \kappa$.

\qed

In the following we investigate this case in another way.
For the Riemannian manifold with constant sectional curvature $\kappa$, $\forall ~ W, X, Y, Z \in T_pM $ we have \cite[P.149]{Lee_Book97}
$$ R(W,X,Y,Z) = \langle W, R(X,Y)Z \rangle = \kappa ( \langle W,X \rangle
 \langle Z,Y \rangle - \langle W,Y \rangle \langle Z,X \rangle ), $$
equivalently,
$$ R(X,Y) Z =  \kappa ( \langle Z,Y \rangle X - \langle Z,X \rangle Y ). $$
Assume that $X$ is the tangent vector of geodesic and $\langle X,X \rangle =1$, and $Y$ is a Jacobi field along the geodesic, normal to $X$, and $\langle Y,Y \rangle =1$. Then setting $Z=X$, we have
\begin{equation}\label{eqn:ConstCurvature1}
R(X,Y) X = \kappa ( \langle X,Y \rangle X - \langle X,X \rangle Y ) = - \kappa Y,
\end{equation}
where we use the facts that $\langle X,X \rangle =1$ and $\langle X,Y \rangle =0$.
Similarly, we can calculate that
\begin{equation}\label{eqn:ConstCurvature2}
R(Y,X) Y = \kappa ( \langle Y,X \rangle Y - \langle Y,Y \rangle X ) = - \kappa X.
\end{equation}
Since $R(X,Y) = -R(Y,X)$, from these two equations \eqref{eqn:ConstCurvature1} and \eqref{eqn:ConstCurvature2}  we can clearly see that $\kappa$ is the M-eigenvalue.

\subsection{The Schwarzschild solution }

The Schwarzschild solution is the first exact solution of Einstein's field equation.
According to Birkhoff's theorem the Schwarzschild solution is the most general spherically symmetric solution of the vacuum Einstein equation.
The exterior Schwarzschild metric is framed in a spherical coordinate system with the body's centre located at the origin, plus the time coordinate.
In Schwarzschild coordinates, with signature $(-1, 1, 1, 1)$, the line element for the Schwarzschild metric has the form
$$ds^2 = - \left( 1-\frac{2GM}{c^2 r} \right) c^2 dt^2 + \left( 1-\frac{2GM}{c^2 r} \right)^{-1} dr^2 + r^2 d\theta^2 + r^2\sin^2 \theta d\phi^2.$$

Let $r_s = \frac{2 GM} {c^2}$ be the Schwarzschild radius, and define
 $f(r) = 1 - \frac{2 GM} {c^2 r} = 1 - \frac{r_s}{r}$.
The nonzero components of the Riemann tensor $R^i_{jkl}$ are given as follows.
$$R^0_{101} = \frac{r_s}{ r^3 f(r) },   \quad
R^0_{220} = R^1_{221} =\frac{r_s}{2 r}, \quad
R^0_{330} = R^1_{331} = \frac{1}{2} R^2_{323} = \frac{r_s }{2 r}\sin^2 \theta , $$
$$R^1_{001} = - g^{11} g_{00} R^0_{101} = \frac{ r_s }{r^3}c^2 f(r),    \quad
R^2_{020} = - g^{22} g_{00} R^0_{220} = \frac{ r_s }{2 r^3}c^2 f(r), $$
$$R^3_{030} = - g^{33} g_{00} R^0_{330} = \frac{ r_s }{2 r^3}c^2 f(r),  \quad
R^2_{112} = g^{22} g_{11} R^1_{221} = \frac{r_s}{2 r^3 f(r)} , $$
$$R^3_{113} = g^{33} g_{11} R^1_{331} = \frac{r_s}{2 r^3 f(r)},         \quad
R^3_{232} = g^{33} g_{22} R^2_{323} = \frac{r_s }{ r}, $$
$$R^0_{110} = -R^0_{101}, ~\qquad  R^0_{202} = -R^0_{220}, ~\qquad R^0_{303} = -R^0_{330}, $$
$$R^1_{010} = -R^1_{001}, ~\qquad  R^1_{212} = -R^1_{221}, ~\qquad R^1_{313} = -R^1_{331}, $$
$$R^2_{002} = -R^2_{020}, ~\qquad  R^2_{121} = -R^2_{112}, ~\qquad R^2_{332} = -R^2_{323}, $$
$$R^3_{003} = -R^3_{030}, ~\qquad  R^3_{131} = -R^3_{113}, ~\qquad R^3_{223} = -R^3_{232}. $$
%

\begin{proposition}
For the Schwarzschild solution given above, $\zeta =  - \sqrt{\frac{K_1}{48}} 
$ is the M-eigenvalue,
where
$K_1 = R_{abcd} R^{abcd} = 12 \left( \frac{r_s}{r^3} \right)^2 = \frac{48 G^2 M^2 }{c^4 r^6} $ is the Kretschman curvature invariant.
\end{proposition}
{\bf Proof}.
Let $t=x^0, r=x^1, \theta = x^2, \phi=x^3$. For the Schwarzschild space, the M-eigenproblem reads
\begin{eqnarray*}
R^0_{101}y^1 (x^0y^1-x^1y^0) +R^0_{202}y^2 (x^0y^2-x^2y^0) +R^0_{303}y^3(x^0y^3-x^3y^0) =\zeta x^0, \\
R^1_{001}y^0 (x^0y^1-x^1y^0) +R^1_{212}y^2 (x^1 y^2-x^2y^1) +R^1_{313} y^3(x^1 y^3-x^3y^1) = \zeta x^1, \\
R^2_{002} y^0 (x^0 y^2-x^2y^0) +R^2_{112} y^1 (x^1 y^2-x^2y^1) +R^2_{323} y^3 (x^2 y^3-x^3y^2) = \zeta x^2, \\
R^3_{003} y^0 (x^0 y^3-x^3y^0) +R^3_{113} y^1 (x^1 y^3-x^3y^1) +R^3_{223} y^2 (x^2 y^3-x^3y^2) = \zeta x^3, \\
R^0_{101}x^1 (y^0x^1-y^1x^0) +R^0_{202}x^2 (y^0x^2-y^2x^0) +R^0_{303}x^3(y^0x^3-y^3x^0) =\zeta y^0, \\
R^1_{001}x^0 (y^0x^1-y^1x^0) +R^1_{212}x^2 (y^1 x^2-y^2x^1) +R^1_{313} x^3(y^1 x^3-y^3x^1) = \zeta y^1, \\
R^2_{002} x^0 (y^0 x^2-y^2x^0) +R^2_{112} x^1 (y^1 x^2-y^2x^1) +R^2_{323} x^3 (y^2 x^3-y^3x^2) = \zeta y^2, \\
R^3_{003} x^0 (y^0 x^3-y^3x^0) +R^3_{113} x^1 (y^1 x^3-y^3x^1) +R^3_{223} x^2 (y^2 x^3-y^3x^2) = \zeta y^3, \\
    g_{00} x^0 x^0 + g_{11} x^1 x^1 + g_{22} x^2 x^2 + g_{33} x^3 x^3 =  1 , \\
    g_{00} y^0 y^0 + g_{11} y^1 y^1 + g_{22} y^2 y^2 + g_{33} y^3 y^3 =  1 .
\end{eqnarray*}

Define  $A=\frac{ GM f(r)}{r^3} =\frac{ r_s }{2 r^3}c^2 f(r)$, $B = \frac{ GM}{ c^2 r^3 f(r) } = \frac{ r_s }{ 2 r^3 f(r) }  $, $C=\frac{GM }{c^2 r} =\frac{r_s }{2 r}$, $D=\frac{GM \sin^2 \theta}{c^2 r} =\frac{r_s }{2 r} \sin^2 \theta  $. The nonzero components are expressed by
\begin{eqnarray*}
 R^0_{101} =-R^0_{110}= 2B, \qquad R^0_{220}=-R^0_{202}= C, \qquad R^0_{330}=-R^0_{303} = D, \\
 R^1_{001}=-R^1_{010} = 2A, \qquad R^1_{221}=-R^1_{212}= C, \qquad R^1_{331}=-R^1_{313} = D ,\\
 R^2_{020} =-R^2_{002}= A,  \qquad R^2_{112}=-R^2_{121}= B, \qquad R^2_{323} =-R^2_{332}= 2D,\\
 R^3_{030}=-R^3_{003} = A,  \qquad R^3_{113}=-R^3_{131}= B, \qquad R^3_{232}=-R^3_{223} = 2C .
\end{eqnarray*}

We have the following system of polynomial equations with nine variables $(x^0, x^1, x^2, x^3,  y^0, y^1, y^2, y^3, \zeta)$.
\begin{eqnarray*}
2B y^1 (x^0y^1-x^1y^0)  - C y^2 (x^0y^2-x^2y^0)  - D y^3(x^0y^3-x^3y^0) =\zeta x^0, \\
2A y^0 (x^0y^1-x^1y^0)  - C y^2 (x^1 y^2-x^2y^1) - D y^3(x^1 y^3-x^3y^1) = \zeta x^1, \\
-A y^0 (x^0 y^2-x^2y^0) + B y^1 (x^1 y^2-x^2y^1) +2D y^3 (x^2 y^3-x^3y^2) = \zeta x^2, \\
-A y^0 (x^0 y^3-x^3y^0) + B y^1 (x^1 y^3-x^3y^1) -2C y^2 (x^2 y^3-x^3y^2) = \zeta x^3, \\
-2B x^1 (x^0y^1-x^1y^0)  + C x^2 (x^0y^2-x^2y^0)  + D x^3(x^0y^3-x^3y^0) =\zeta y^0, \\
-2A x^0 (x^0y^1-x^1y^0)  + C x^2 (x^1 y^2-x^2y^1) + D x^3(x^1 y^3-x^3y^1) = \zeta y^1, \\
 A x^0 (x^0y^2-x^2y^0) - B x^1 (x^1 y^2-x^2y^1) - 2D x^3 (x^2 y^3-x^3y^2) = \zeta y^2, \\
 A x^0 (x^0y^3-x^3y^0) - B x^1 (x^1 y^3-x^3y^1) + 2C x^2 (x^2 y^3-x^3y^2) = \zeta y^3.
\end{eqnarray*}

Let $S^{ij} = x^i y^j - x^j y^i$.
From the first four equalities, we have
$$
\begin{bmatrix}
0 & 2S^{01} & S^{20} & S^{30} \\ 2S^{01} & 0 & S^{21} & S^{31} \\
S^{20} & S^{12} & 0 & 2S^{23} \\ S^{30} & S^{13} & 2S^{32} & 0 \\
\end{bmatrix}
\begin{bmatrix}
A & 0 & 0 & 0 \\ 0 & B & 0 & 0 \\ 0 & 0 & C & 0 \\ 0 & 0 & 0 & D \\
\end{bmatrix}
\begin{bmatrix}
y^0 \\ y^1 \\ y^2 \\ y^3
\end{bmatrix} = \zeta
\begin{bmatrix}
x^0 \\ x^1 \\ x^2 \\ x^3
\end{bmatrix}
$$
For simplicity, we denote it by
$$ T \Lambda y = \zeta x , $$
where $\Lambda :=  \text{diag}(A,B,C,D)  =  \frac{r_s} {2 r^3} \text{diag}(c^2 f(r), f(r)^{-1}, r^2, r^2 \sin^2 \theta) $, and $T$ can be defined correspondingly.

Similarly, from the last four equalities, we have
$$ T \Lambda x = - \zeta y .$$

Rewriting these two equations in a compact form, we have
$$
\begin{bmatrix} 0 & T \Lambda \\ - T \Lambda & 0 \end{bmatrix}
\begin{pmatrix} x \\ y \end{pmatrix} = \zeta
\begin{pmatrix} x \\ y \end{pmatrix} .
$$

Note that
$$ \det
\begin{pmatrix} \zeta & - T \Lambda \\   T \Lambda & \zeta \end{pmatrix}
= \det \left(  \zeta^2 I + (T \Lambda)^2 \right)
 . $$
The nonzero solution requires that  $\det \left( T \Lambda \pm i \zeta I   \right) =0$. Direct calculation yields that
\begin{eqnarray*}
\det \left( T  \pm i \zeta \Lambda^{-1}   \right) =   \det(T) +
[ \zeta^2 + A B (2S^{01})^2 + A C (S^{20})^2 + A D (S^{30})^2  \\
 - C D (2S^{23})^2 - B D (S^{13})^2 - B C (S^{21})^2   ] \zeta^2  \det(\Lambda^{-1})  .
\end{eqnarray*}
A sufficient condition for the existence of nonzero solution is that the two terms in the right hand size of above formula are zeros.

For the first term, we can verify that
$$ \det(T) = - ( 2 S^{01} \cdot 2 S^{23} + S^{20} S^{13} + S^{30} S^{21} )^2 .$$
Hence, $ \det(T) = 0$ is equivalent to 
$4 S^{01} S^{23} + S^{20} S^{13} + S^{30} S^{21} =0 $.
We can denote this by
$ (2 S^{01}, S^{20}, S^{30} ) \perp (2 S^{23}, S^{13}, S^{21} )$ for an easy-to-remember form. Substituting the express for $S^{ij}$, we have $3 S^{01} S^{23} = 0$, i.e., $x^0 y^1 = x^1 y^0$ or $x^2 y^3=x^3y^2$. Meanwhile, $ S^{20} S^{31} - S^{30} S^{21} = 0$.

The choice that $x=y$ gives zero M-eigenvalue. We choose
$$ x^0 = y^0 , \quad x^1 = y^1, \quad x^2 = -y^2, \quad x^3 = -y^3. $$
Under such settings, we have $S^{01}=0=S^{23}$, $S^{20}= 2 x^0 x^2$, $S^{30}= 2 x^0 x^3$, $S^{13}= -2 x^1 x^3$, and $S^{21}= 2 x^1 x^2$.
For this special choice, the first four equations and the next four equations are the same, and can be further reduced to two equations:
\begin{eqnarray*}
A (x^0)^2 - B (x^1)^2 &=&  \zeta /2 ,\\
C (x^2)^2 + D (x^3)^2 &=& -\zeta /2 ,
\end{eqnarray*}
which can be obtained, for example, from the 2nd and the 3rd equations in the system of polynomial equations.
This indicates that
$$A(x^0)^2 - B(x^1)^2 = - C(x^2)^2 - D(x^3)^2  .$$

By the way, we can check that the second term in $\det \left( T  \pm i \zeta \Lambda^{-1}   \right)$ reduces to
$$ \{  \zeta^2 + 4 [ A(x^0)^2 - B(x^1)^2 ] [ C (x^2)^2 + D (x^3)^2 ] \} \zeta^2  \det(\Lambda^{-1}), $$
and it vanishes naturally.

Since the metric $$(g_{ij}) = \text{diag}(-c^2 f(r), f(r)^{-1}, r^2, r^2 \sin^2 \theta) = \frac{2 r^3}{r_s} \text{diag}(-A, B, C, D),  $$
and the constraint
$$\langle x, x \rangle = g_{ij} x^i x^j = \frac{2r^3}{r_s} [- A(x^0)^2 + B(x^1)^2 + C(x^2)^2 + D(x^3)^2] =  1, $$
we have
$$ 2 [ A(x^0)^2 - B(x^1)^2 ] = \zeta = - \frac{r_s}{2 r^3}  = - \frac{GM}{c^2 r^3} = - \sqrt{ \frac{K_1}{48} } .$$

\qed

We find that the M-eigenvalue relates to the Kretschmann invariant, which  is the simplest invariant product involving the Riemann
curvature tensor, and used most often to identify essential singularities in a spacetime
geometry.

The Schwarzschild metric (1915), the Reissner--Nordstr\"{o}m metric (1916, 1918), the Kerr metric (1963), and the Kerr--Newman metric (1965) are four related solutions.
The analysis above can be extended to the Reissner--Nordstr\"{o}m solution, a static solution 
corresponding to the gravitational field of a charged, non-rotating, spherically symmetric body. Its line element reads
$$ds^2 = - \left( 1-\frac{r_s}{r} + \frac{r_Q^2}{r^2}  \right) c^2 dt^2 + \left( 1-\frac{r_s}{r} + \frac{r_Q^2}{r^2} \right)^{-1} dr^2  + r^2 ( d\theta^2 + \sin^2 \theta d\phi^2 ) , $$
where $r_Q$ is a characteristic length scale given by
$r_Q^2 = \frac{G Q^2}{4\pi \epsilon_0 c^4}$ with $Q$ being the charge. The Kerr solution is more interesting, but the corresponding M-eigenproblem is much more difficult.

\section{Conclusion}

Starting from two kinds of eigenproblems of the elasticity tensor, we investigate two corresponding eigenproblems of the Riemann curvature tensor. One is classical, related to the Petrov classification. The other one is the M-eigenvalue problem, which is the extension of the M-eigenvalue problem of the elasticity tensor. In a way similar to the elasticity tensor case, from the optimization problem on the associated function, we introduce the M-eigenproblem of the Riemann tensor. The M-eigenvalues of the Riemann curvature tensor always exist and are real.   They are invariants of the Riemann curvature tensor.   The associated function of the Riemann curvature tensor is always positive at a point if and only if the M-eigenvalues of the Riemann curvature tensor are all positive at that point.    These show that the M-eigenvalues are some intrinsic scalars of the Riemann curvature tensor.  We further examine several typical cases such as the 2D case, the 3D case, the constant curvature case and the Schwarzschild solution. But we can just obtain a few of the M-eigenvalues, and cannot calculate all the M-eigenvalues by solving a system of polynomial equations.   Actually, the M-eigenvalues and their corresponding eigenvectors are real solutions of a system of polynomial equations.   Thus, it is hard to know the number of M-eigenvalues \cite{Gelfand_Book1994,Qi_JSC05}. But we believe that the M-eigenvalues are related to those important curvature invariants, and would show us tremendous mathematical and physical information.

\section*{Acknowledgements}
The authors would like to thank the two anonymous referees and Xiaokai He for the helpful comments to improve the presentation of this work.

\section*{Appendix A. }

Consider the case where $z^{ij} = x^i y^j + y^i x^j$ in \eqref{eqn:ElasticityEigentensor}. The problem \eqref{eqn:ElasticityEigentensor} reads
\begin{equation*}
E_{ijkl} ( x^k y^l + y^k x^l ) = \zeta ~ ( x_i y_j + y_i x_j ) .
\end{equation*}
It is easy to check that
\begin{eqnarray*}
E_{ijkl} y^j x^k y^l  + E_{ijkl}  y^j y^k x^l  = \zeta ~ ( y^j x_i y_j + y^j y_i x_j  ) .
\end{eqnarray*}
Using $E_{ijkl}  y^j y^k x^l = E_{ikjl}  y^k y^j x^l = E_{iklj} y^k x^l y^j = E_{ijkl} y^j x^k y^l$, we obtain
\begin{eqnarray} \label{eqn:ElasticityEig0}
2 E_{ijkl} y^j x^k y^l  = \zeta ~ ( x_i  +  \langle x, y \rangle y_i  ).
\end{eqnarray}

The problem \eqref{eqn:ElasticityEigentensor}   is equivalent to $E_{ijkl} z^{ij} = \zeta z_{kl}$. Under the same condition, it can be rewritten as
\begin{equation*}
E_{ijkl} ( x^i y^j + y^i x^j ) = \zeta ~ ( x_k y_l + y_k x_l ) .
\end{equation*}
And we can check that
\begin{eqnarray*}
E_{ijkl} x^i y^j x^k + E_{ijkl}  y^i x^j x^k  = \zeta ~ ( x_k y_l x^k + x_l y_k x^k  ) .
\end{eqnarray*}
Using $E_{ijkl}  y^i x^j x^k = E_{ikjl}  y^i x^k x^j = E_{kijl} x^k y^i x^j = E_{ijkl} x^i y^j x^k$, we then have
\begin{eqnarray}\label{eqn:ElasticityEig1}
2 E_{ijkl} x^i y^j x^k  = \zeta ~ ( y_l  +  \langle x, y \rangle x_l  ).
\end{eqnarray}

Both \eqref{eqn:ElasticityEig0} and \eqref{eqn:ElasticityEig1} yield
\begin{equation*}
\zeta = \frac{ 2 E_{ijkl} x^i y^j x^k y^l } {   1  +  \langle x, y \rangle^2  } .
\end{equation*}

When $\langle x,y \rangle = 0$,
\eqref{eqn:ElasticityEig0} and \eqref{eqn:ElasticityEig1}  reduce to the M-eigenvalue problem \eqref{eqn:MeigDef_x&y} with $(\frac{1}{2}\zeta, x, y)$ being the M-eigentriple.
What's more,
$\zeta = 2 E_{ijkl} x^i y^j x^k y^l = 2 \theta $.

\section*{Appendix B. }

Using the symmetry of Riemann tensor $R_{ijkl}$ and the anti-symmetry of eigentensor $x^{ij}$, we can rewrite the eigenproblem \eqref{eqn:RCurvatureEigentensor} as follows.
\begin{eqnarray*}
\begin{bmatrix}
R_{1010} & R_{1020} & R_{1030} & R_{1023} & R_{1031} & R_{1012} \\
R_{2010} & R_{2020} & R_{2030} & R_{2023} & R_{2031} & R_{2012} \\
R_{3010} & R_{3020} & R_{3030} & R_{3023} & R_{3031} & R_{3012} \\
R_{2310} & R_{2320} & R_{2330} & R_{2323} & R_{2331} & R_{2312} \\
R_{3110} & R_{3120} & R_{3130} & R_{3123} & R_{3131} & R_{3112} \\
R_{1210} & R_{1220} & R_{1230} & R_{1223} & R_{1231} & R_{1212}
\end{bmatrix}
\begin{pmatrix}
x^{10} \\ x^{20} \\ x^{30} \\ x^{23} \\ x^{31} \\ x^{12}
\end{pmatrix} =
\frac{1}{2} \zeta
\begin{bmatrix}
G_{11}   & G_{12} \\
G_{21}   & G_{22}
\end{bmatrix}
\begin{pmatrix}
x^{10} \\ x^{20} \\ x^{30} \\ x^{23} \\ x^{31} \\ x^{12}
\end{pmatrix}.
\end{eqnarray*}

We express it in a compact matrix form,  compressing each pair of indices into one index by the following standard mapping for tensor indices.
$$ \begin{matrix}ij & = \\ \Downarrow & \\ A  & = \end{matrix} \quad \begin{matrix} 10 &20 &30 &23 &31 &12 \\ \Downarrow &\Downarrow &\Downarrow &\Downarrow &\Downarrow &\Downarrow &\\1&2&3&4&5&6\end{matrix} $$
Denoting the basis indices by capital letters, we rewrite the above eigenproblem as \eqref{eqn:RCurvatureEig_inMatrix}:
$R_{AB} x^B = \frac{1}{2} \zeta G_{AB} x^B $,
where $(G_{AB})$ is defined by the following matrix blocks of size 3-by-3:
$$G_{11} = \begin{bmatrix}
g_{00}g_{11}-g_{01}^2   & g_{00}g_{12}-g_{02}g_{10} & g_{00}g_{13}-g_{03}g_{10} \\
                        & g_{00}g_{22}-g_{02}^2     & g_{00}g_{23}-g_{03}g_{20}  \\
sym.                    &                           & g_{00}g_{33}-g_{03}^2
\end{bmatrix},
$$
$$G_{22} = \begin{bmatrix}
g_{22}g_{33}-g_{23}^2   & g_{23}g_{31}-g_{21}g_{33} & g_{21}g_{32}-g_{22}g_{31} \\
                        & g_{11}g_{33}-g_{13}^2     & g_{12}g_{31}-g_{11}g_{32}  \\
sym.                    &                           & g_{11}g_{22}-g_{12}^2
\end{bmatrix},
$$
$$G_{12} = \begin{bmatrix}
g_{03}g_{12}-g_{02}g_{13} & g_{01}g_{13}-g_{03}g_{11} & g_{02}g_{11}-g_{01}g_{12} \\
g_{03}g_{22}-g_{02}g_{23} & g_{01}g_{23}-g_{03}g_{21} & g_{02}g_{21}-g_{01}g_{22} \\
g_{03}g_{32}-g_{02}g_{33} & g_{01}g_{33}-g_{03}g_{31} & g_{02}g_{31}-g_{01}g_{32}
\end{bmatrix} ,
$$
and $ G_{21}= (G_{12})^T$.
For the case where $(g_{ij}) = \text{diag} (-1, 1,1,1)$, we have $G = \text{diag} (-1,-1,-1,1,1,1)$.

From the Einstein field equation in vacuum
$$R_{ij} = \kappa g_{ij} ,$$
 and $R_{ab} = R^h_{ahb} = g^{hl} R_{a l b h}$, we have
$$\Sigma_k s_k R_{ikjk} = \kappa g_{ij}, \quad s_k = \pm 1  .$$
Since $R$ is symmetric, we suppose that $ ( R_{AB} ) = \begin{bmatrix} M & N \\ N^T & W \end{bmatrix} $, where $ M$ and $ W$ are symmetric. According to the (0,1), (0,2) and (0,3) components of the equation, we can derive that $N$ itself is symmetric.
Using the Bianchi identity $R_{0123} + R_{0312} + R_{0231} = 0$, we have $\text{tr}(N) = 0 $.
From the (1,2), (1,3) and (2,3) components, we obtain that $M = - W$ except the diagonals. Further more, using the (0,0), (1,1), (2,2) and (3,3) components, we can verify that $ M = - W$ and $ \text{tr} (M) = - \kappa $.
In all, the 6-by-6 matrix $ ( R_{AB} )$ can be expressed as
$$ ( R_{AB} ) = \begin{bmatrix} M & N \\ N & -M \end{bmatrix} , $$
where $\text{tr} (N) =0$, $\text{tr} (M) = - \kappa $, $N$ and $M$ are both symmetric.

\end{document}